\documentclass[10pt,leqno]{article}

\usepackage{amsmath,amssymb,amsthm,mathrsfs,dsfont}

\usepackage[margin=3cm]{geometry} %宽版时用，配合双倍行距

\usepackage{titlesec,hyperref}

\usepackage{color}
\usepackage{fancyhdr}
\titleformat{\subsection}{\it}{\thesubsection.\enspace}{1pt}{}

\newtheorem{theo}{Theorem}[section]
\newtheorem{lemm}[theo]{Lemma}
\newtheorem{defi}[theo]{Definition}

\newtheorem{rema}[theo]{Remark}
\numberwithin{equation}{section}

\allowdisplaybreaks %允许公式分页显示

\newcommand\ep{{\varepsilon}} %可定义一个简单符号来代替很长而常用的命令。

\begin{document}
\title{Blow-up phenomenon, ill-posedness and peakon solutions for the periodic Euler-Poincar\'e equations
\hspace{-4mm}
}

\author{Wei $\mbox{Luo}^1$\footnote{E-mail:  luowei23@mail2.sysu.edu.cn} \quad and\quad
 Zhaoyang $\mbox{Yin}^{1,2}$\footnote{E-mail: mcsyzy@mail.sysu.edu.cn}\\
 $^1\mbox{Department}$ of Mathematics,
Sun Yat-sen University,\\ Guangzhou, 510275, China\\
$^2\mbox{Faculty}$ of Information Technology,\\ Macau University of Science and Technology, Macau, China}
\date{}
\maketitle
\hrule

\begin{abstract}
In this paper we mainly investigate the initial value problem of the periodic Euler-Poincar\'e equations. We first present a new blow-up result to the system for a special class of smooth initial data by using the rotational invariant properties of the system. Then, we prove that the periodic Euler-Poincar\'e equations is ill-posed in critical Besov spaces by a contradiction argument. Finally, we verify the system possesses a class of peakon solutions in the sense of distributions.\\

\vspace*{5pt}

\vspace*{5pt}
\end{abstract}

\vspace*{10pt}

%\phantomsection
%\addcontentsline{toc}{section}{\contentsname}l
%添加目录到书签
\tableofcontents
\section{Introduction}
  In this paper we consider the initial value problem for the following periodic Euler-Poinc\'are equations:
  \begin{align}\label{1}
\left\{
\begin{array}{ll}
m_{t}+u\cdot\nabla m+ (\nabla u)^{T} m+ m~(div~{u})=0,\\[1ex]
m|_{t=0}(x)=m_{0}(x),\\[1ex]
m(t,x)=m(t,x+1),\\[1ex]
\end{array}
\right.
\end{align}
or in components,
\begin{equation}
\left\{
\begin{array}{ll}
m^i_t+\sum^{d}_{j=1}u^j\partial_j m^i+\sum^d_{j=1}m^j\partial_i u^j+m_i\sum^d_{j=1}\partial_j u^j=0,\quad i=1,2,...,d, \\[1ex]
m^i|_{t=0}(x)=m^i_0(x).\\[1ex] m(t,x)=m(t,x+1).
\end{array}
\right.
\end{equation}
Here $u$ is the velocity and $m=(1-\Delta)u$ denotes the momentum, $d$ is the spatial dimension. The Euler-Poincar\'e equations were first studied in
\cite{Holm1998-Adv,Holm1998-Phy,Holm2003-SIAM} as a higher dimensional Camassa-Holm system for modeling and analyzing the nonlinear shallow water waves. Moreover, it can be viewed as the geodesic flow on $\mathcal{D}^1(M)$ with the right-invariant metric \cite{Holm1998-Adv,Holm1998-Phy}:
\[\langle u,v\rangle=\int_{M}uv+u_xv_x dx,\]
where $M$ is a Riemannian manifold (In general, $M$ is $\mathbb{R}^d$ or $\mathbb{T}^d$) . If $d=1$ the Euler-Poincar\'e equations are reduced to the following famous Camassa-Holm equation \cite{Camassa1993}:
\begin{align}
m_t+2m_xu+mu_x=0, ~~m=u-u_{xx}.
\end{align}
The Camassa-Holm equation was derived as a model for shallow water waves \cite{Camassa1993, Constantin2009}. It has been investigated extensively because of its great physical significance in the past two decades. The CH equation has a bi-Hamiltonian structure \cite{Constantin1997,Fokas1981} and is completely integrable \cite{Camassa1993,Constantin2001}. The solitary wave solutions of the CH equation were considered in \cite{Camassa1993,Camassa1994}, where the authors showed that the CH equation possesses peakon solutions of the form $Ce^{-|x-Ct|}$. It is worth mentioning that the peakons are solitons and their shape is alike that of the travelling water waves of greatest height, arising as solutions to the free-boundary problem for incompressible Euler equations over a flat bed (these being the governing equations for water waves),
cf. the discussions in \cite{Constantin2006,Constantin2007,Constantin2011,Toland1996}. Constantin and Strauss verified that the peakon solutions of the CH equation are orbitally stable in \cite{Constantin.S2000}. \\
 $~~~~~~$The local well-posedness for the CH equation was studied in \cite{Constantin1998,Constantin.1998,Danchin2001,Blanco2001,Li-Yin-2016-JDE}. Concretely, for initial profiles $u_0\in H^s(\mathbb{R})$ with $s>\frac{3}{2}$, it was shown in \cite{Constantin1998,Constantin.1998,Blanco2001} that the CH equation has a unique solution in $C([0,T);H^s(\mathbb{R}))$. Moveover, the local well-posedness for the CH equation in Besov spaces $C([0,T);B^s_{p,r}(\mathbb{R}))$ with $s>\max(\frac{3}{2},1+\frac{1}{p})$ was proved in \cite{Danchin2001,Li-Yin-2016-JDE}. The global existence of strong solutions were established in \cite{Constantin2000,Constantin1998,Constantin.1998} under some sign conditions and it was shown in \cite{Constantin2000,Constantin1998,Constantin.1998,Constantin-1998} that the solutions will blow up in finite time when the slope of initial data was bounded by a negative quantity. The global weak solutions for the CH equation were studied in \cite{Constantin.M} and \cite{Xin2000}. The global conservative and dissipative solutions of CH equation were presented in \cite{Bressan2007,Holden2007} and  \cite{Bressan.2007}, respectively.\\

Recently, the local well-posedness in $H^s(\mathbb{R}^d)$ with $s>1+\frac{d}{2}$ for the Euler-Poincar\'e equations was established in \cite{Chae-2012-CMP}. Moreover, the authors obtained a blow-up criteria, zero $\alpha$ limit and the Liouville type theorem. In \cite{Li-2013-ARMA}, the authors constructed a special class of solutions ($m=\nabla\phi$) to the Euler-Poincar\'e equations which will blow up in finite time. The local well-posedness for the high-dimension CH equations in Besov spaces $C([0,T);B^s_{p,r}(\mathbb{R}^d))$ with $s>\max(\frac{3}{2},1+\frac{d}{p})$ or $s=1+\frac{d}{p},~r=1$ was proved in \cite{Yan-2015-DCDS}.

The authors in \cite{Li-2013-ARMA} apply the symmetrical structure to reduce the higher dimension system into a one dimension equation, and then prove that this kind of special solutions to the Euler-Poincar\'e equations will blow up in finite time. In this paper we construct a new special class of strong solutions to the periodic Euler-Poincar\'e equations (\ref{1}) which will also blow up in finite time. Our approach and the obtained result are quite different from the recent result in \cite{Li-2013-ARMA}. The main idea is that we use the rotational invariant properties of the system (\ref{1}). If the initial data is constant along some vector filed, and then the solution has the same property. This observation leads us to obtain a new blow-up result.

Recently, the authors in \cite{Guo-Liu-Yin-2018-JDE} prove ill-posedness of the Camassa-Holm type equations in the critical spaces. To our best knowledge, there is no any ill-posedness results for the Euler-Poincar\'e equations. Inspired by the works of \cite{Guo-Liu-Yin-2018-JDE,Linares} about the Camassa-Holm equation and the Burgers equation, we use the contradiction argument to prove the system (\ref{1}) is ill-posed in critical Besov spaces.

The remainder of the paper is organized as follows. In Section 2 we introduce some preliminaries which will be used in sequel. In Section 3 we prove the blow-up phenomenon of the system (\ref{1}) with $d=2$. In Section 4 we generalize the obtained result from $d=2$ to $d\geq 3$. In Section 5, we prove that the system (\ref{1}) is ill-posed in critical Besov spaces. Section 6 is devoted to verifying that the system (\ref{1}) possesses a class of peakon solutions.

\section{Preliminaries}
In this section, we recall some previous results for the Euler-Poincar\'e equations.

In \cite{Chae-2012-CMP}, the authors write $(\nabla u)^{T}m$ in a tensor form:
\begin{align*}
\sum^d_{j=1}m^j\partial_i u^j &=\sum^d_{j=1}u^j\partial_iu^j-\sum^{d}_{j,k=1}\partial^2_k u^j\partial_iu^j=\frac{1}{2}\partial_i|u|^2-\sum^{d}_{j,k=1}\partial_k(\partial_iu^j\partial_ku^j)+\sum^{d}_{j,k=1}\partial_k\partial_iu^j\partial_ku^j\\
&=\sum^d_{j=1}\partial_j(\frac{1}{2}\delta_{ij}|u|^2-\partial_i u\cdot\partial_j u+\frac{1}{2}\delta_{ij}|\nabla u|^2).
\end{align*}
Denote the tensor $T^{ij}=m^iu^j+\frac{1}{2}\delta_{ij}|u|^2-\partial_i u\cdot\partial_j u+\frac{1}{2}\delta_{ij}|\nabla u|^2$. Then the system (\ref{1}) becomes
\[\partial_t m^i+\sum^d_{j=1}\partial_jT^{ij}=0.\]

Let us recall the local well-posedness result and the blow-up criteria for the Euler-Poincar\'e equations.
\begin{lemm}\cite{Chae-2012-CMP}
Let $u_0\in H^k(\mathbb{R}^d)$ with $k>\frac{d}{2}+3$. Then, there exists $T=T(\|u_0\|_{H^k(\mathbb{R}^d)})$ such that the Euler-Poincar\'e equations have a unique solution $u\in C([0,T];H^k(\mathbb{R}^d)$.
\end{lemm}

\begin{lemm}\cite{Chae-2012-CMP}
Let $u_0\in H^k(\mathbb{R}^d)$ with $k>\frac{d}{2}+3$. Suppose that $T^*$ is the lifespan of the solution $u$ to the Euler-Poincar\'e equations. Then the solution blows up in finite time if and only if
\[\int^{T^*}_0\|S(t)\|_{\dot{B}^0_{\infty,\infty}(\mathbb{R}^d)}dt=\infty,\]
where $S=(S^{ij})$ with $S^{ij}=\frac{1}{2}\partial_iu^j+\frac{1}{2}\partial_ju^i.$
\end{lemm}

\begin{rema}
The embedding relation $L^\infty(\mathbb{R}^d)\hookrightarrow BMO(\mathbb{R}^d)\hookrightarrow \dot{B}^0_{\infty,\infty}(\mathbb{R}^d) $ implies that
\[\|S(t)\|_{\dot{B}^0_{\infty,\infty}(\mathbb{R}^d)}\leq \|S\|_{L^\infty(\mathbb{R}^d)}\leq \|\nabla u\|_{L^\infty(\mathbb{R}^d)}.\]
Therefore we obtain the following criteria:
\[\limsup_{t\rightarrow T^*}\|u\|_{H^k(\mathbb{R}^d)}=\infty \quad \text{if and only if} \quad \int^{T^*}_0\|\nabla u\|_{L^\infty(\mathbb{R}^d)}dt=\infty.\]
\end{rema}

For the periodic case, by virtue of the transport equation theory one can obtain the similar results as follows:
\begin{lemm}\label{L2.4}
Let $u_0\in H^k(\mathbb{T}^d)$ with $k>\frac{d}{2}+3$. Then, there exists $T=T(\|u_0\|_{H^k(\mathbb{T}^d)})$ such that the system (\ref{1}) has a unique solution $u\in C([0,T];H^k(\mathbb{T}^d)$.
\end{lemm}

\begin{lemm}
Let $u_0\in H^k(\mathbb{T}^d)$ with $k>\frac{d}{2}+3$. Suppose that $T^*$ is the lifespan of the solution $u$ to (\ref{1}). Then the solution blows up in finite time if and only if
\[\int^{T^*}_0\|\nabla u\|_{L^\infty}dt=\infty.\]
\end{lemm}
Since the regularity index $k>\frac{d}{2}+3$, it follows that the space $H^{k-2}(\mathbb{T}^d)$ is a Banach algebra, which ensures the local existence and uniqueness of the solution to a transport equation. Since the proof is similar to that of \cite{Chae-2012-CMP}, we omit the details here.

Now, we introduce a crucial lemma which will be used to prove our main result.
\begin{lemm}\label{L6}
Let $u_0\in H^k(\mathbb{T}^d)$ with $k>\frac{d}{2}+3$ and let $u$ be the corresponding local solution to (\ref{1}). If $\frac{\partial u^i_0}{\partial n}=0,~i=1,2,...d$, then $\frac{\partial u^i}{\partial n}=0,~~t\in[0,T]$, where $n=(n^1,n^2,...,n^d)\in \mathbb{R}^d$ is a constant vector field.
\begin{proof}
Since $\frac{\partial f}{\partial n}=\sum^d_{i=1}n^i\partial_if$, it follows that
\[\frac{\partial (fg)}{\partial n}=\frac{\partial f}{\partial n}g+f\frac{\partial g}{\partial n}.\]
Applying $\frac{\partial }{\partial n}$ to both sides of (\ref{1}) in components, we have
\begin{align}
\partial_nm^i_t+\sum^d_j(u^j\partial^2_{jn}m^i+\partial_nu^j\partial_jm^i+\partial^2_{in}u^jm^j+\partial_iu^j\partial_nm^j+\partial_nm^i\partial_ju^j+m^i\partial^2_{nj}u^j)=0.
\end{align}
Multiplying by $\partial_n m^i$ both sides of the above inequality and integrating over $\mathbb{T}^d$, we obtain
\begin{align}\label{2.2}
&\frac{d}{dt}\int_{\mathbb{T}^d}|\partial_n m^i|^2dx=-\sum^d_{j}\int_{\mathbb{T}^d}(\partial_nu^j\partial_jm^i+\partial^2_{in}u^jm^j+\partial_iu^j\partial_nm^j+m^i\partial^2_{nj}u^j)\partial_n m^idx\\
\nonumber&\leq (\|\nabla m\|_{L^\infty}\|\partial_n u\|_{L^2}+2\|m\|_{L^\infty}\|\partial_n\nabla u\|_{L^2}+\|\nabla u\|_{L^\infty}\|\partial_n m\|_{L^2})\|\partial_n m\|_{L^2}.
\end{align}
Since $m=(1-\Delta )u$, it follows that
\[\|\partial_n u\|_{L^2}\leq \|\partial_n m\|_{L^2},~~\|\partial_n\nabla u\|_{L^2}\leq \|\partial_n m\|_{L^2}.\]
Plugging the above inequality into (\ref{2.2}) yields that
\begin{align}
\frac{d}{dt}\|\partial_nm\|^2_{L^2}\leq (\|\nabla m\|_{L^\infty}+2\|m\|_{L^\infty}+\|\nabla u\|_{L^\infty})\|\partial_nm\|^2_{L^2}\leq C\|u\|_{H^k}\|\partial_nm\|^2_{L^2},
\end{align}
where we use the embedding relation $H^k\hookrightarrow W^{3,\infty}$ with $k>\frac{d}{2}+1$. Taking advantage of Gronwall's inequality and the fact that $\partial_nm_0=0$, we see that
\[\|\partial_n m\|_{L^2}=0.\]
\end{proof}
\end{lemm}
Similar to the whole space case \cite{Chae-2012-CMP, Li-2013-ARMA}, one can easily get the following lemma.
\begin{lemm}
The system \ref{1} has the following conservation law:
\begin{align}
H=\int_{\mathbb{T}^d}|u|^2+|\nabla u|^2 dx.
\end{align}
\end{lemm}

For a function $f$ on $\mathbb{T}^d$, we define its Fourier transform denoted by $\widehat{f}(\xi)$ as
\[\widehat{f}(\xi)=\frac{1}{(2\pi)^d}\int_{\mathbb{T}^d}f(x)e^{-i\xi x}dx, \quad \xi\in \mathbb{Z}^d.\]
Let $\eta:\mathbb{R} \rightarrow [0,1]$ be an even, smooth, non-negative and radially decreasing function which is supported in a ball $\{\xi:|\xi|\leq \frac{8}{5}\}$ and $\eta\equiv 1$ for $|\xi|\leq \frac{5}{4}$.
Let $\varphi(\xi)\doteq\eta(\xi)-\eta(\frac{\xi}{2})$, and define the Littlewood-Paley operators $\Delta_j$ by \[\widehat{\Delta_j f}(\xi)\doteq\varphi(2^{-j}\xi)\widehat{f}(\xi), ~~\widehat{\Delta_{-1} f}(\xi)\doteq\eta(\xi)\widehat{f}(\xi),~~ S_j f\doteq\sum_{j'<j}\Delta_{j'} f \quad \forall j\in\mathbb{Z},\]
We can then define the Besov space $B^s_{p,r}(\mathbb{T}^d)$ with norm
\[\|f\|_{B^s_{p,r}(\mathbb{T}^d)}=\|2^{js}\|\Delta_j f\|_{L^p(\mathbb{T}^d)}\|_{l^r}.\]

{\bf Notation.} Since all function spaces in the following sections are over $\mathbb{T}^d$, for simplicity, we drop $\mathbb{T}^d$ in the notation of function spaces if there is no ambiguity.
\section{Blow up: $d=2$}
In order to explain our main idea to obtain the blow-up result, we first consider a simple case $d=2$. Our main theorem in this section can be stated as follows.
\begin{theo}
 Let $u_0\in H^k(\mathbb{T}^2)$ with $k>4$ and let $u$ be the corresponding local solution to (\ref{1}). Suppose that $\frac{\partial u^i_0}{\partial n}=0,~~i=1,2$,  where $n=(\cos\theta, \sin\theta)$, $\theta$ is a constant.  If $\sin^2\theta  \partial_x u^1_0-\sin\theta\cos\theta\partial_y u^1_0-\sin\theta\cos\theta\partial_xu^2_0+\cos^2\theta \partial_yu^2_0<-\sqrt{2}\|u_0\|_{H^1}$ for some $(x_0,y_0)\in\mathbb{T}^2$. Then the solution $u$ blows up in finite time.
 \begin{proof}
 Let $n^{\bot}=(-\sin\theta,\cos\theta)$ be the vertical vector of $n$. Define that
 \begin{align}
 u^n=u\cdot n=\cos\theta u^1+\sin\theta u^2,~~\quad u^{n^{\bot}}=u \cdot n^{\bot}=-\sin\theta u^1+\cos\theta u^2,
 \end{align}
 and then
\begin{align}
u^1=\cos\theta u^n-\sin\theta u^{n^{\bot}},~~\quad u^2=\sin\theta u^{n}+\cos\theta u^{n^{\bot}}.
\end{align}
 By directly calculating, we see that
 \begin{align}
 div~u=\partial_xu^1+\partial_yu^2=\partial_x(\cos\theta u^n-\sin\theta u^{n^{\bot}})+\partial_y (\sin\theta u^{n}+\cos\theta u^{n^{\bot}})=\frac{\partial u^n}{\partial n}+\frac{\partial u^{n^{\bot}}}{\partial n^{\bot}},
 \end{align}
 \begin{align}
 u\cdot \nabla=u^1\partial_x+u^2\partial_y=(\cos\theta u^n-\sin\theta u^{n^{\bot}})\partial_x+(\sin\theta u^{n}+\cos\theta u^{n^{\bot}})\partial_y=u^n\frac{\partial}{\partial n}+u^{n^{\bot}}\frac{\partial }{\partial n^{\bot}}.
 \end{align}
Multiplying by $n^{\bot}$ both sides of (\ref{1}), we obtain
\begin{align}\label{3.5}
\frac{\partial m^{n^{\bot}}}{\partial t}+u^n\frac{\partial m^{n^{\bot}}}{\partial n }+u^{n^{\bot}}\frac{\partial m^{n^{\bot}}}{\partial n^{\bot}}+m^{n^{\bot}}(\frac{\partial u^n}{\partial n}+\frac{\partial u^{n^{\bot}}}{\partial n^{\bot}})+(\nabla u)^{T}m\cdot n^{\bot}=0.
\end{align}
Since $\frac{\partial u^i_0}{\partial_n}=0$, by virtue of Lemma (\ref{L6}), we have $\frac{\partial u^i}{\partial_n}=0$ and
\begin{align}\label{3.6}
\frac{\partial m^{n^{\bot}}}{\partial t}+u^{n^{\bot}}\frac{\partial m^{n^{\bot}}}{\partial n^{\bot}}+m^{n^{\bot}}\frac{\partial u^{n^{\bot}}}{\partial n^{\bot}}+(\nabla u)^{T}m\cdot n^{\bot}=0.
\end{align}
Next we deal with the term $(\nabla u)^{T}m\cdot n^{\bot}$ in components.
\begin{align}\label{3.7}
&(\nabla u)^{T}m\cdot n^{\bot}=-\sin\theta(m^1\partial_xu^1+m^2\partial_xu^2)+\cos\theta(m^1\partial_yu^1+m^2\partial_yu^2)\\
\nonumber&=-\sin\theta\{\frac{1}{2}\partial_x[(u^1)^2-(\partial_xu^1)^2+(\partial_yu^1)^2+(u^2)^2-(\partial_xu^2)^2+(\partial_yu^2)^2]-\partial_y(\partial_yu^1\partial_xu^1+\partial_yu^2\partial_xu^2)\}\\
\nonumber&+\cos\theta\{\frac{1}{2}\partial_y[(u^1)^2-(\partial_yu^1)^2+(\partial_xu^1)^2+(u^2)^2-(\partial_yu^2)^2+(\partial_xu^2)^2]-\partial_x(\partial_yu^1\partial_xu^1+\partial_yu^2\partial_xu^2)\}\\
\nonumber&=\frac{1}{2}\frac{\partial[(u^1)^2+(u^2)^2]}{\partial{n^{\bot}}}+(\frac{\sin\theta}{2}\partial_x+\frac{\cos\theta}{2}\partial_y)[(\partial_xu^1)^2+(\partial_xu^2)^2-(\partial_yu^1)^2+(\partial_yu^2)^2]\\
\nonumber&+(\sin\theta\partial_y-\cos\theta\partial_x)(\partial_yu^1\partial_xu^1+\partial_yu^2\partial_xu^2)\\
\nonumber&=\frac{1}{2}\frac{\partial (u^n)^2}{\partial{n^{\bot}}}+\frac{1}{2}\frac{\partial (u^{n^\bot})^2}{\partial{n^{\bot}}}+(\frac{\sin\theta}{2}\partial_x+\frac{\cos\theta}{2}\partial_y)[(\partial_xu^1)^2+(\partial_xu^2)^2-(\partial_yu^1)^2+(\partial_yu^2)^2]\\
\nonumber&+(\sin\theta\partial_y-\cos\theta\partial_x)(\partial_yu^1\partial_xu^1+\partial_yu^2\partial_xu^2),
\end{align}
where we use the fact that $(u^1)^2+(u^2)^2=(u^n)^2+(u^{n^{\bot}})^2$. Since $\partial_nu^1=\partial_nu^2=0$, it follows that
\begin{align}\label{3.8}
\sin\theta\partial_x+\cos\theta\partial_y=(\cos^2\theta-\sin^2\theta)\frac{\partial}{\partial n^\bot},~~\sin\theta\partial_y-\cos\theta\partial_x=2\sin\theta\cos\theta\frac{\partial}{\partial n^\bot},
\end{align}
\begin{align}\label{3.9}
(\partial_xu^1)^2+(\partial_xu^2)^2=(-\sin\theta\frac{\partial u^1}{\partial n^{\bot}})^2+(-\sin\theta\frac{\partial u^2}{\partial n^{\bot}})^2=\sin^2\theta[(\frac{\partial u^n}{\partial n^{\bot}})^2+(\frac{\partial u^{n^{\bot}}}{\partial n^{\bot}})^2],\\
\nonumber (\partial_yu^1)^2+(\partial_yu^2)^2=(\cos\theta\frac{\partial u^1}{\partial n^{\bot}})^2+(\cos\theta\frac{\partial u^2}{\partial n^{\bot}})^2=\cos^2\theta[(\frac{\partial u^n}{\partial n^{\bot}})^2+(\frac{\partial u^{n^{\bot}}}{\partial n^{\bot}})^2],\\
\nonumber \partial_yu^1\partial_xu^1+\partial_yu^2\partial_xu^2=-\cos\theta\sin\theta[(\frac{\partial u^1}{\partial n^{\bot}})^2+(\frac{\partial u^2}{\partial n^{\bot}})^2]=-\cos\theta\sin\theta[(\frac{\partial u^n}{\partial n^{\bot}})^2+(\frac{\partial u^{n^{\bot}}}{\partial n^{\bot}})^2].
\end{align}
Plugging (\ref{3.8}) and (\ref{3.9}) into (\ref{3.7}) yields that
\begin{align}\label{3.10}
(\nabla u)^{T}m\cdot n^{\bot}=\frac{1}{2}\frac{\partial (u^n)^2}{\partial{n^{\bot}}}+\frac{1}{2}\frac{\partial (u^{n^{\bot}})^2}{\partial{n^{\bot}}}-\frac{1}{2}\frac{\partial}{\partial n^{\bot}}[(\frac{\partial u^n}{\partial n^{\bot}})^2+(\frac{\partial u^{n^{\bot}}}{\partial n^{\bot}})^2].
\end{align}
For simplicity, we denote that $v=\partial_{n^{\bot}} u^{n^{\bot}}$, $w=\partial_{n^{\bot}} u^{n}$. Using (\ref{3.6}) and (\ref{3.10}) we get
\begin{align}\label{3.11}
\frac{\partial m^{n^{\bot}}}{\partial t}+u^{n^{\bot}}\frac{\partial m^{n^{\bot}}}{\partial n^{\bot}}+m^{n^{\bot}}\frac{\partial u^{n^{\bot}}}{\partial n^{\bot}}+\frac{1}{2}\frac{\partial (u^n)^2}{\partial{n^{\bot}}}+\frac{1}{2}\frac{\partial (u^{n^{\bot}})^2}{\partial{n^{\bot}}}
-\frac{1}{2}\frac{\partial}{\partial n^{\bot}}(v^2+w^2)=0.
\end{align}
Since $m=(1-\Delta)u=(1-\partial^2_x-\partial^2_y)u=(1-\partial^2_{ n^{\bot}})u$, it follows that
\begin{align}\label{3.12}
(1-\partial^2_{n^{\bot}})(\frac{\partial u^{n^{\bot}}}{\partial t}+u^{n^{\bot}}\frac{\partial u^{n^{\bot}}}{\partial n^{\bot}})=\frac{\partial m^{n^{\bot}}}{\partial t}+u^{n^{\bot}}\frac{\partial m^{n^{\bot}}}{\partial n^{\bot}}+3m^{n^{\bot}}\frac{\partial u^{n^{\bot}}}{\partial n^{\bot}}-3u^{n^{\bot}}\frac{\partial u^{n^{\bot}}}{\partial n^{\bot}}.
\end{align}
Combining (\ref{3.12}) with (\ref{3.11}) yields that
\begin{align}\label{3.13}
\frac{\partial u^{n^{\bot}}}{\partial t}+u^{n^{\bot}}\frac{\partial u^{n^{\bot}}}{\partial n^{\bot}}=-\partial_{n^{\bot}}(1-\partial^2_{n^{\bot}})^{-1}\bigg[(u^{n^{\bot}})^2+\frac{1}{2}(u^n)^2+\frac{1}{2}v^2-\frac{1}{2}w^2\bigg].
\end{align}
Applying $\frac{\partial }{\partial n^{\bot}}$ on both sides of the above inequality,  we have
\begin{align}
\partial_tv+u^{n^\bot}\partial_{n^\bot}v&=-v^2-\partial^2_{n^{\bot}}(1-\partial^2_{n^{\bot}})^{-1}f=-v^2+f-(1-\partial^2_{n^{\bot}})^{-1}f\\
\nonumber&=-\frac{1}{2}(v^2+w^2)+(u^{n^{\bot}})^2+\frac{1}{2}(u^n)^2-(1-\partial^2_{n^{\bot}})^{-1}f,
\end{align}
where $f=(u^{n^{\bot}})^2+\frac{1}{2}(u^n)^2
+\frac{1}{2}v^2-\frac{1}{2}w^2$.

Define the characteristics $\Phi(t,x,y)$ as the solution of
\begin{align}\label{O}
\left\{
\begin{array}{ll}
\frac{d\Phi(t,x,y)}{dt}=u(t,\Phi(t,x,y)),\\[1ex]
\Phi(0,x,y)=(x,y).
\end{array}
\right.
\end{align}
It is easy to check that
\[\frac{d}{dt}v(t,\Phi(t,x,y))=v_t+u\nabla v=\partial_tv+u^{n^\bot}\partial_{n^\bot}v,~~~\]
which leads to
\begin{align}\label{3.16}
\frac{d}{dt}v(t,\Phi(t,x,y))=-\frac{1}{2}(v^2+w^2)(t,\Phi(t,x,y))+(u^{n^\bot})^2+\frac{1}{2}(u^n)^2-(1-\partial^2_{n^{\bot}})^{-1}f.
\end{align}
Since $\|u\|_{H^1}=\|u_0\|_{H^1}$, it follows that
\begin{align}
-(1-\partial^2_{n^{\bot}})^{-1}f\leq \frac{1}{2}(1-\partial^2_{n^{\bot}})^{-1}w^2\leq \frac{1}{2}\|u_0\|^2_{H^1}.
\end{align}
By Sobolev's embedding, we see that $\|u^n\|^2_{L^\infty}+\|u^{n^\bot}\|^2_{L^\infty}\leq \frac{1}{2}(\|\frac{\partial u^n}{\partial n^{\bot}}\|^2_{L^2}+\|\frac{\partial u^{n^\bot}}{\partial n^{\bot}}\|^2_{L^2})\leq \frac{1}{2}\|u\|^2_{H^1}=\frac{1}{2}\|u_0\|^2_{H^1}$.
Using the fact that $-\frac{1}{2}w^2\leq 0$ and the above estimate, we have
\begin{align}
\frac{d g(t)}{dt}\leq -\frac{1}{2}g^2(t)+\|u_0\|^2_{H^1}\leq  -\frac{1}{2}(g(t)+\sqrt{2}\|u_0\|_{H^1})(g(t)-\sqrt{2}\|u_0\|_{H^1}),
\end{align}
where $g(t)=v(t,\Phi(t,x_0,y_0))$. The assumption on the initial datum and the definition of $v$ guarantee that $g(0)<-\sqrt{2}\|u_0\|_{H^1}$. By virtue of the continuity argument, we see that $g(t)<-\sqrt{2}\|u_0\|_{H^1}$ for all $t\in[0,T^*)$. From the above inequality, we obtain that
\begin{align}
-\frac{d}{dt}(\frac{1}{g(t)-\sqrt{2}\|u_0\|_{H^1}})\leq -(1+\frac{2\sqrt{2}\|u_0\|_{H^1}}{g(t)-\sqrt{2}\|u_0\|_{H^1}}),
\end{align}
which leads to
\begin{align}
\frac{g(0)+\sqrt{2} \|u_0\|_{H^1}}{g(0)-\sqrt{2}\|u_0\|_{H^1}}e^{\sqrt{2}\|u_0\|_{H^1}t}-1\leq \frac{2\sqrt{2}\|u_0\|_{H^1}}{g(t)-\sqrt{2}\|u_0\|_{H^1}}\leq 0.
\end{align}
Since $0<\frac{g(0)+\sqrt{2}\|u_0\|_{H^1}}{g(0)-\sqrt{2}\|u_0\|_{H^1}}<1$, then there exists
\[0<T\leq \frac{\sqrt{2}}{2 \|u_0\|_{H^1}}\ln\frac{g(0)-\sqrt{2}\|u_0\|_{H^1}}{g(0)+\sqrt{2} \|u_0\|_{H^1}},\]
such that $\lim_{t\rightarrow T}g(t)=-\infty$. By the definition of $g(t)$, we see that $|g(t)|\leq \|v\|_{L^\infty}\leq \|\frac{\partial u^{n^\bot}}{\partial n^\bot}\|_{L^\infty}\leq \|\nabla u\|_{L^\infty}$, which implies that the solution blows up in finite time.
 \end{proof}
\end{theo}

\begin{rema}
If a non-zero function $f$ satisfies that $\partial_nf=0$, then the $H^1$ norm of $f$ in $\mathbb{R}^2$ is infinity i.e. $\|f\|_{H^1(\mathbb{R}^2)}=\infty$. Thus, the argument in the above theorem only holds true in the periodic case. For the whole space, one should use the other function space to deal with this problem.
\end{rema}

\section{Blow up: $d\geq 3$}
Now we turn our attention to the higher dimension ($d\geq 3$) case. Let's first introduce some notations. Let
\[A=(A_1,A_2,...A_d)=\{(a^{ij})\},\]
where $A_j$ is the column vector of the matrix of $A$. Denote that
\[\overline{u}=u^A=Au, \quad \overline{\nabla}f=(A \nabla)f,\]
or in components
\[\overline{u}^i=\sum_{j}a^{ij}u^j, \quad \overline{\partial_i}=\sum_{j}a^{ij}\partial_j.\]

Our main result in this section can be stated as follows.
\begin{theo}\label{th2}
Assume that $A$ is an orthogonal matrix. Let $u_0\in H^k(\mathbb{T}^d)$ with $k>\frac{d}{2}+3$ and let $u$ be the corresponding local solution to (\ref{1}). Suppose that $\overline{\partial_i}u_0=0,~~i=2,3,...,d$.  If there exists some $x_0\in \mathbb{T}^d$ such that $\overline{\partial_1}\overline{u}^1_0<-\sqrt2\|u_0\|_{H^1}$. Then the solution $u$ blows up in finite time.
\begin{proof}
Since $A$ is an orthogonal matrix, it follows that
\begin{align}
div~u= \sum^d_{i=1}\partial_i u^i=\sum^d_{i,j=1}a^{ji}\overline{\partial_j}u^i=\sum^d_{j=1}\overline{\partial_j}\overline{u}^j=\overline{div} ~\overline{u}, \\
u\cdot \nabla=\sum^d_{i=1}u^i\partial_i =\sum^d_{i,j=1}a^{ji}\overline{u^j}\partial_i=\sum^d_{j=1}\overline{u}^j\overline{\partial_j}=\overline{u}\cdot\overline{\nabla},
\end{align}
Applying $A$ to both sides of (\ref{1}) and using the above equality, we see that
\begin{align}\label{4.3}
\overline{m}_{t}+\overline{u}\cdot\overline{\nabla} \overline{m}+ A\cdot[(\nabla u)^{T} m]+ \overline{m}~(\overline{div}~\overline{u})=0.
\end{align}
According to Lemma \ref{L6} and the assumption $\overline{\partial_i}u_0=0,~i=2,3,...,d$ , we get $\overline{\partial_i}u=0,~i=2,3,...,d$. From (\ref{4.3}), we see that
\begin{align}\label{4.4}
\overline{m}_{t}+\overline{u}^1\cdot\overline{\partial_1} \overline{m}+ A\cdot[(\nabla u)^{T} m]+ \overline{m}~(\overline{\partial_1}\overline{u}^1)=0.
\end{align}
Now we consider the term $A\cdot[(\nabla u)^{T} m]$ in components.  By directly calculating, we deduce that
\begin{align}\label{4.5}
( A\cdot[(\nabla u)^{T} m])^i&= \sum^d_{j,k=1}a^{ij}m^k\partial_j u^k=\sum^d_{j,k=1}a^{ij}\partial_k(\frac{1}{2}\delta_{jk}|u|^2-\partial_j u\cdot\partial_k u+\frac{1}{2}\delta_{jk}|\nabla u|^2)\\
\nonumber&=\frac{1}{2}\sum^d_{j=1}a^{ij}\partial_j(|u|^2+|\nabla u|^2)-\sum^d_{j,k=1}a^{ij}\partial_k(\partial_j u\cdot\partial_k u)\\
\nonumber&=\frac{1}{2}\overline{\partial_i}(|u|^2+|\nabla u|^2)-\sum^d_{j,k=1}a^{ij}\partial_k(\partial_j u\cdot\partial_k u).
\end{align}
Since $A$ is an orthogonal matrix, it follows that
\begin{align}\label{4.6}
|u|^2=u^T\cdot u=u^TA^{T}Au=(Au)^T(Au)=|\overline{u}|^2, \\
\nonumber|\nabla u|^2=(\nabla u)^T\cdot\nabla u=(\nabla u)^TA^{T}A\nabla u=|\overline{\nabla}\overline{u}|^2=|\overline{\partial_1}\overline{u}|^2.
\end{align}
By directly calculating, we see that
\begin{align}\label{4.7}
\partial_k(\partial_j u\cdot\partial_k u)=\sum^d_{l=1}a^{lk}\overline{\partial_l}[\sum^d_{i=1}(\sum^d_{m=1}a^{jm}\overline{\partial_m}\overline{u}^i)(\sum^d_{n=1}a^{kn}\overline{\partial_n}\overline{u}^i)]=a^{1k}a^{j1}a^{k1}\overline{\partial_1}(|\overline{\partial_1}\overline{u}|^2).
\end{align}
Combining (\ref{4.6}) and (\ref{4.7}) with (\ref{4.5}) yields that
\begin{align}\label{4.8}
( A\cdot[(\nabla u)^{T} m])^1
=\frac{1}{2}\overline{\partial_1}(|\overline{u}|^2+|\overline{\partial_1} \overline{u}|^2)-\sum^d_{j,k=1}a^{1j}a^{1k}a^{j1}a^{k1}\overline{\partial_1}(|\overline{\partial_1}\overline{u}|^2)=\frac{1}{2}\overline{\partial_1}(|\overline{u}|^2-|\overline{\partial_1} \overline{u}|^2),
\end{align}
where we use the fact that $\sum_{k}a^{1k}a^{k1}=1$. From (\ref{4.4}) and (\ref{4.8}) we deduce that
\begin{align}\label{4.9}
\overline{m}^1_{t}+\overline{u}^1\cdot\overline{\partial_1} \overline{m}^1+\frac{1}{2}\overline{\partial_1}(|\overline{u}|^2-|\overline{\partial_1} \overline{u}|^2) + \overline{m}^1~(\overline{\partial_1}\overline{u}^1)=0.
\end{align}
Since $\overline{m}^1=(1-\overline{\partial^2_1})\overline{u}^1$, it follows that
\begin{align}\label{4.10}
(1-\overline{\partial^2_1})(\overline{u}^1_t+\overline{u}^1\overline{\partial_1}\overline{u}^1)=\overline{m}^1_t+\overline{u}^1\overline{\partial_1}\overline{m}^1+3\overline{m}^1~(\overline{\partial_1}\overline{u}^1)-3 \overline{u}^1~(\overline{\partial_1}\overline{u}^1).
\end{align}
Plugging (\ref{4.10}) into (\ref{4.9}) yields that
\begin{align}
\overline{u}^1_t+\overline{u}^1\overline{\partial_1}\overline{u}^1=-\overline{\partial_1}(1-\overline{\partial^2_1})^{-1}[\frac{1}{2}(\overline{u}^1)^2+\frac{1}{2}|\overline{u}|^2-\frac{1}{2}|\overline{\partial_1} \overline{u}|^2+(\overline{\partial_1} \overline{u})^2].
\end{align}
Applying $\overline{\partial_1}$ to both sides of the above inequality, we have
\begin{align}
\overline{\partial_1}\overline{u}^1_t+\overline{u}^1\overline{\partial_1}\overline{u}^1&=-(\overline{\partial_1}\overline{u}^1)^2-\overline{\partial^2_1}(1-\overline{\partial^2_1})^{-1}f=-(\overline{\partial_1}\overline{u}^1)^2+f-(1-\overline{\partial^2_1})^{-1}f\\
\nonumber&=-\frac{1}{2}|\overline{\partial_1} \overline{u}|^2+\frac{1}{2}(\overline{u}^1)^2+\frac{1}{2}|\overline{u}|^2-(1-\overline{\partial^2_1})^{-1}f,
\end{align}
where $f=\frac{1}{2}(\overline{u}^1)^2+\frac{1}{2}|\overline{u}|^2-\frac{1}{2}|\overline{\partial_1} \overline{u}|^2+(\overline{\partial_1} \overline{u})^2$.
Define the characteristics $\Phi(t,x)\in \mathbb{R}^d$ as the solution of
\begin{align}\label{O1}
\left\{
\begin{array}{ll}
\frac{d\Phi(t,x)}{dt}=u(t,\Phi(t,x)),\\[1ex]
\Phi(0,x)=x.
\end{array}
\right.
\end{align}
Then we deduce that
\[\frac{d}{dt}v(t,\Phi(t,x))=v_t+u\nabla v=\partial_tv+\overline{u}^1\overline{\partial_1}v,~~~\]
which leads to
\begin{align}
\frac{d}{dt}(\overline{\partial_1}\overline{u}^1\circ\Phi)=-\frac{1}{2}|\overline{\partial_1} \overline{u}\circ\Phi|^2+\frac{1}{2}(\overline{u}^1\circ\Phi)^2+\frac{1}{2}|\overline{u}\circ\Phi|^2-(1-\overline{\partial^2_1})^{-1}f\circ\Phi.
\end{align}
By virtue of the conservation law $\|u\|_{H^1}=\|u_0\|_{H^1}$, we obtain
\begin{align}
-(1-\overline{\partial^2_1})^{-1}f\leq \frac{1}{2}(1-\overline{\partial^2_1})^{-1}|\overline{\partial_1} \overline{u}|^2\leq \frac{1}{2}\|u_0\|^2_{H^1}.
\end{align}
Taking advantage of Sobolev's embedding, we see that $\|\overline{u}\|^2_{L^\infty}\leq \frac{1}{2}\|\overline{\partial_1}\overline{u}\|^2_{L^2}\leq \frac{1}{2}\|u\|^2_{H^1}=\frac{1}{2}\|u_0\|^2_{H^1}$.
Using the fact that $-\frac{1}{2}|\overline{\partial_1}\overline{u}|^2\leq -\frac{1}{2}(\overline{\partial_1}\overline{u}^1)^2$ and the above estimate, we have
\begin{align}\label{4.16}
\frac{d h(t)}{dt}\leq -\frac{1}{2}h^2(t)+\|u_0\|^2_{H^1}\leq  -\frac{1}{2}(h(t)+\sqrt{2}\|u_0\|_{H^1})(h(t)-\sqrt{2}\|u_0\|_{H^1}),
\end{align}
where $h(t)=\overline{\partial_1}\overline{u}^1\circ\Phi(t,x_0)$. Note that $h(0)<-\sqrt{2}\|u_0\|_{H^1}$. The continuity argument ensures that $h(t)<-\sqrt{2}\|u_0\|_{H^1}$ for all $t\in[0,T^*)$. By the same token, we can solve the inequality (\ref{4.16}) and deduce that
\begin{align}
\frac{h(0)+\sqrt{2} \|u_0\|_{H^1}}{h(0)-\sqrt{2}\|u_0\|_{H^1}}e^{\sqrt{2}\|u_0\|_{H^1}t}-1\leq \frac{2\sqrt{2}\|u_0\|_{H^1}}{h(t)-\sqrt{2}\|u_0\|_{H^1}}\leq 0.
\end{align}
Thanks to $0<\frac{h(0)+\sqrt{2}\|u_0\|_{H^1}}{h(0)-\sqrt{2}\|u_0\|_{H^1}}<1$, we can find
\begin{align}\label{4.18}
0<T\leq \frac{\sqrt{2}}{2 \|u_0\|_{H^1}}\ln\frac{h(0)-\sqrt{2}\|u_0\|_{H^1}}{h(0)+\sqrt{2} \|u_0\|_{H^1}},
\end{align}
such that $\lim_{t\rightarrow T}h(t)=-\infty$. Since $|h(t)|\leq \|\overline{\partial_1}\overline{u}^1\|_{L^\infty}\leq \|\nabla u\|_{L^\infty}$, it follows that the solution blows up in $T$.
\end{proof}
\end{theo}

\section{Ill-posedness in critical Besov spaces}
In this section, we are going to prove the norm inflation of the system (\ref{1}) in  $B^{1+\frac{d}{p}}_{p,r}$ with $d\leq p\leq +\infty$ and $1<r\leq +\infty$. Our main result can be sated as follows:
\begin{theo}\label{th3}
Let $d\leq p\leq +\infty$ and $1<r\leq +\infty$. For any $\ep>0$, there exists $u_0\in H^\infty(\mathbb{T}^d)$, such that the following holds:\\
$(1)$ $\|u_0\|_{B^{1+\frac{d}{p}}_{p,r}(\mathbb{T}^d)}\leq \ep$;\\
$(2)$ There exists a unique solution $u\in C([0,T);H^\infty(\mathbb{T}^d))$ of the system (\ref{1}) with maximal $T<\ep$;\\
$(3)$  $\limsup_{t\rightarrow T^-} \|u\|_{B^{1+\frac{d}{p}}_{p,r}(\mathbb{T}^d)}\geq \limsup_{t\rightarrow T^-} \|u\|_{B^{1}_{\infty,\infty}(\mathbb{T}^d)}=\infty$.
\end{theo}

\begin{rema}
$(i)$ In the case $d=2$, taking $p=r=2$, we see that our theorem implies the ill-posedness in the critical Sobolev space $H^{2}$. \\
$(ii)$ If $d\geq 3$, we don't know whether the system (\ref{1}) is well-posed or ill-posed in the critical Sobolev space $H^{1+\frac{d}{2}}$.
\end{rema}

In order to prove Theorem \ref{th3}, we need some useful lemmas.
\begin{lemm}\label{L5.3}\cite{Guo-Liu-Yin-2018-JDE}
Let $T>0$. Assume that $A(t)\in C^1[0,T), A(t)>0$ and there exists a constant $B$ such that
\[ \frac{d}{dt}A(t)\leq BA(t)\ln(2+A(t)), \quad \forall t\in[0,T).\]
Then we have
\[A(t)\leq (2+A(0))^{e^{Bt}}, \quad \forall t \in[0,T).\]
\end{lemm}

\begin{lemm}\label{L5.4}\cite{Guo-Liu-Yin-2018-JDE}
Assume that $u\in H^{1+\frac{d}{2}+\ep}$ with $\ep>0$. We have
\[\|\nabla u\|_{L^\infty} \leq \frac{C}{\ep}\|u\|_{B^1_{\infty,\infty}}\ln(2+\frac{\|u\|_{H^{1+\frac{d}{2}+\ep}}}{\|u\|_{B^1_{\infty,\infty}}})\]
\end{lemm}

{\bf Proof of Theorem \ref{th3}}: From now on, we turn our attention to prove Theorem \ref{th3}.
Fix $d\leq p\leq \infty, 1< r\leq  \infty$. We define the periodic function $h(x_1)$ as follows:
\[W(x_1)=-\sum_{k\geq 2} \frac{1}{k^{1+\frac{d}{p}}\ln^{\frac{2}{1+r}} k}\sin(kx_1).\]
Directly calculating, we have
\[\Delta_j W(x_1)= -\sum_{2^j<k<2^{j+2}}\varphi(\frac{k}{2^j}\vec{e}_1)\frac{1}{k^{1+\frac{d}{p}}\ln^{\frac{2}{1+r}} k}\sin(kx_1),\]
where $\vec{e}_1=(1,0,0,...,0)$. On the other hand,
\[\|\sin(kx_1)\|_{L^p(\mathbb{T}^d)}=(\int_{T^d}|\sin(kx_1)|^pdx)^{\frac{1}{p}}=(2\pi)^{\frac{d-1}{p}}(\frac{1}{k}\int^{2k\pi}_0|\sin(y)|^pdy)^{\frac{1}{p}}=C_{p,d},\]
where $C_{p,d}$ is a constant independent on $k$, which implies that $\|\Delta_j W(x_1)\|_{L^p(\mathbb{T}^d)}\sim 2^{-j(1+\frac{d}{p})}j^{-\frac{2}{1+r}}$, and thus
\[\|W(x_1)\|_{B^{1+\frac{d}{p}}_{p,q}(\mathbb{T}^d)}\sim \|\frac{1}{j^{\frac{2}{1+r}}}\|_{l^q}.\]
From this we see that $h \in B^{1+\frac{d}{p}}_{p,r}$ but $h\notin B^{1+\frac{d}{p}}_{p,1}$. Since $\frac{d}{p}\leq 1$ and $r>1$, it follows that
\[W'(0)=-\sum_{k\geq 2} \frac{1}{k^{\frac{d}{p}}\ln^{\frac{2}{1+r}} k} =-\infty.\]

For any $\ep>0$, let
\[u^1_{0,\ep}=\displaystyle\frac{\ep S_N W(x_1)}{\|W\|_{B^{1+\frac{d}{p}}_{p,r}}}, ~~~u^2_0=u^3_0=...=u^n_0=\ep f(x_1), \]
where $N$ is large enough such that $\partial_1u^1_{0,\ep}(0)\sim \ep^{-10}$ and $f$ is arbitrary smooth periodic function such that $\|f\|_{B^{1+\frac{d}{p}}_{p,r}}<1$.
(For example, one can choose $f=\frac{\sin(x_1)}{\|\sin(x_1)\|_{B^{1+\frac{d}{p}}_{p,r}}}$).
Then $u_{0,\ep}\in H^\infty$ and $\|u_0\|_{B^{1+\frac{d}{p}}_{p,r}}\leq \ep$.
By virtue of Lemma \ref{L2.4}, there exists a $T_\ep$ such that the system (\ref{1}) has a unique solution $u_\ep\in C([0,T_\ep); H^\infty)$.  Taking $A=Id$ and $x_0=0$ in Theorem \ref{th2}, from the estimate (\ref{4.18}), we see that the maximal time
\[T_\ep  \leq \frac{\sqrt{2}}{2 \|u_{0,\ep}\|_{H^1}}\ln\frac{\partial_1u^1_{0,\ep}(0)-\sqrt{2}\|u_{0,\ep}\|_{H^1}}{\partial_1u^1_{0,\ep}(0)+\sqrt{2} \|u_{0,\ep}\|_{H^1}}\leq C\ep^{10}.\]

It suffices to show that
\begin{align}\label{5.1}
\limsup_{t\rightarrow T_\ep} \|u_\ep\|_{B^{1}_{\infty,\infty}(\mathbb{T}^d)}=\infty.
\end{align}

Suppose that (\ref{5.1}) fails, we can find a $M_\ep$ such that
\[\sup_{0<t<T_\ep}\|u_\ep\|_{B^{1}_{\infty,\infty}}\leq M_\ep.\]

Let $k=2+[\frac{d}{2}]>1+\frac{d}{2}$. Taking advantage of the energy estimate and Lemma \ref{L5.4}, we obtain
\begin{align}
\frac{d}{dt}\|u_\ep\|^2_{H^{k}}\leq \|\nabla u\|_{L^\infty}\|u\|^2_{H^k}\leq C_k\|u_\ep\|_{B^{1}_{\infty,\infty}}\|u\|^2_{H^k}(\ln(2+\frac{\|u\|_{H^k})}{\|u_\ep\|_{B^{1}_{\infty,\infty}}})
\leq C_kM_\ep\|u\|^2_{H^k} \ln(2+\|u\|^2_{H^k}).
\end{align}
By virtue of Lemma \ref{L5.3}, we deduce that $\sup_{t\in[0,T_\ep)}\|u_\ep\|_{H^k}<+\infty$ which means that $u_\ep$ will not blow up in $T_\ep$. This contradicts with Theorem \ref{th2}.

\section{Periodic peakon solutions}
In this section we verify that the system (\ref{1}) possesses a special class of periodic peakon solutions as follows:
\begin{align}\label{6.1}
u^i=M\Phi(a\cdot x-Ct)=M\Phi(\sum^d_{j=1}a^jx^j -Ct), i=1,2,...d,
\end{align}
where $\Phi(z)=\frac{\cosh{(\frac{1}{2}-z)}}{\sinh{(\frac{1}{2})}}$ is a periodic function with $z\in[0.1]$, $M$, $C$ are constants and $a=(a^1,a^2,...a^d)$ is a vector such that $C=\frac{\cosh(\frac{1}{2})}{\sinh(\frac{1}{2})}M\sum^d_{j=1} a^j$ and $|a|=1$.

Firstly, we recall the definition for weak solution of the system (\ref{1}).
\begin{defi}\cite{Chae-2012-CMP}
$u\in L^\infty(0,T;H^1(\mathbb{T}^d))$ is a weak solution of the system (\ref{1}) with initial data $u_0\in H^1(\mathbb{T}^d)$ if the following equation holds for all vector field $\phi(t,x)$ such that $\phi\in C^1([0,T);C^\infty(\mathbb{T}^d)$ and $\Phi(T,x)=0$ for all $x\in\mathbb{T}^d$,
\begin{multline}
\int^T_0\int_{\mathbb{T}^d}(u\cdot\phi_t+\nabla u:\nabla \phi_t) dxdt+\int_{\mathbb{T}^d}(u_0\cdot\phi(0,x)+\nabla u_0:\nabla \phi(0,x))dx\\+\int^T_0\int_{\mathbb{T}^d}T^a:\nabla \phi dxdt
+\sum^d_{i,j,k=1}\int^T_0\int_{\mathbb{T}^d}u^j\partial_ku^i\partial_{jk}\phi^i dxdt=0,
\end{multline}
where $T^a=u\otimes u+\nabla u(\nabla u)^T-(\nabla u)^T(\nabla u)+\frac{1}{2}(|u|^2+|\nabla u|^2)Id$ is the symmetric part of the tensor $T$.
\end{defi}
From (\ref{6.1}), we see that
\begin{align}\label{5.3}
\partial_ku^i=-a^kM\frac{\sinh{(\frac{1}{2}-a\cdot x+Ct+l)}}{\sinh{(\frac{1}{2})}}, \quad \text{if}\quad  a\cdot x-Ct\in(l,l+1),\\
\partial^2_ku^i=-(a^k)^2M\frac{\cosh{(\frac{1}{2}-a\cdot x+Ct+l)}}{\sinh{(\frac{1}{2})}}, \quad \text{if}\quad  a\cdot x-Ct\in(l,l+1).
\end{align}
Note that $\frac{\Phi}{2}$ is the Green function of $1-\partial_{zz}$. By virtue of the above equality, we obtain that $(1-\Delta)u^i=2M\delta(a\cdot x-Ct)$ where $\delta$ is the Dirac function.
Taking advantage of integration by parts, we have
\begin{align}\label{5.5}
\int^T_0\int_{\mathbb{T}^d}(u\cdot\phi_t+\nabla u:\nabla \phi_t) dxdt=2\sum^d_{i=1}\int^T_0\int_{\mathbb{T}^d}M\delta(a\cdot x-Ct)\cdot\phi^i_tdxdt,\\
\int_{\mathbb{T}^d}(u_0\cdot\phi(0,x)+\nabla u_0:\nabla \phi(0,x))dx=2\sum^d_{i=1}\int_{\mathbb{T}^d}M\delta(a\cdot x)\cdot\phi^i(0,x)dx.
\end{align}
Using the fact that $\sinh^2z=\cosh^2z-1$ and (\ref{5.3}), we deduce that
\begin{align}
&[u\otimes u+(\nabla u)(\nabla u)^T-(\nabla u)^T(\nabla u)]:\nabla\phi=\sum^{d}_{i,k}[u^iu^k+\sum^d_j(\partial_ju^i\partial_ju^k-\partial_iu^j\partial_ku^j)]\partial_{i}\phi^k\\
\nonumber&=M^2\sum^d_{i,k=1}(2\Phi^2(a\cdot x-Ct)-1-d\Phi^2(a\cdot x-Ct)a^ia^k-da^ia^k)\partial_{i}\phi^k.
\end{align}
Since $\Phi$ is periodic function, it follows that
\begin{align}\label{5.8}
\int^T_0\int_{\mathbb{T}^d}T^a:\nabla \phi dxdt
&=d M^2 \sum^{d}_{i,k=1}\int^T_0\int_{\mathbb{T}^d}\Phi^2(a\cdot x-Ct)(\delta_{ik}-a^ia^k)\partial_{i}\phi^kdxdt\\
\nonumber&+2M^2\sum^{d}_{i,k=1}\int^T_0\int_{\mathbb{T}^d}\Phi^2(a\cdot x-Ct)\partial_{i}\phi^kdxdt.
\end{align}
Taking advantage of integration by parts, we have
\begin{multline}\label{5.9}
\sum^{d}_{i,k=1}\int^T_0\int_{\mathbb{T}^d}\Phi^2(a\cdot x-Ct)(\delta_{ik}-a^ia^k)\partial_{i}\phi^kdxdt\\
=-\sum^{d}_{i,k=1}\int^T_0\int_{\mathbb{T}^d}(a^i\delta_{ik}-(a^i)^2a^k)\Phi(a\cdot x-Ct)\Phi'(a\cdot x-Ct)\phi^kdxdt=0.
\end{multline}
\begin{multline}\label{5.10}
\sum^d_{k=1}\int^T_0\int_{\mathbb{T}^d}u^j\partial_ku^i\partial_{jk}\phi^i dxdt=-\sum^d_{k=1}\int^T_0\int_{\mathbb{T}^d}(u^j\partial^2_ku^i+\partial_ku^j\partial_ku^i)(\partial_j\phi^i)dxdt\\
-M^2\int^T_0\int_{\mathbb{T}^d}[\Phi(a\cdot x-Ct)\Phi''(a\cdot x-Ct)\partial_j\phi^i+\Phi^2(a\cdot x-Ct)\partial_{j}\phi^i]dxdt.
\end{multline}
Combining (\ref{5.8}), (\ref{5.9}) and (\ref{5.10}),  we see that
\begin{multline}\label{5.11}
\int^T_0\int_{\mathbb{T}^d}T^a:\nabla \phi dxdt+\sum^d_{i,j,k=1}\int^T_0\int_{\mathbb{T}^d}u^j\partial_ku^i\partial_{jk}\phi^i dxdt\\
=2M^2\sum^{d}_{i,k=1}\int^T_0\int_{\mathbb{T}^d}\Phi(a\cdot x-Ct)\delta(a\cdot x-Ct)\partial_i\phi^kdxdt.
\end{multline}

Let $y=a\cdot x$ and $dx=dydz$. By virtue of the properties of the Dirac function, we deduce that
\begin{align}\label{5.12}
&\int^T_0\int_{\mathbb{T}^d}\delta(a\cdot x-Ct)\cdot\phi^i_tdxdt=\int^T_0\int_{\Sigma}\phi^i_t(t,Ct,z)dz,\\
&\int_{\mathbb{T}^d}\delta(a\cdot x)\cdot\phi^i(0,x)dx=\int_{\Sigma}\phi^i(0,0,z)dz,\\
&\int^T_0\int_{\mathbb{T}^d}\Phi(a\cdot x-Ct)\delta(a\cdot x-Ct)\partial_i\phi^kdxdt=\frac{\cosh(\frac{1}{2})}{\sinh(\frac{1}{2})}\int^T_0\int_{\Sigma}\partial_k\phi^i(t,Ct,z)dz.
\end{align}
Note that $\frac{\cosh(\frac{1}{2})}{\sinh(\frac{1}{2})}M(\sum^d_{j}a^j)=C$. Since \[\frac{d}{dt}(\phi^i(t,Ct,z))=\phi^i_t(t,Ct,z)+C\partial_y\phi^i(t,Ct,z),\] it follows that
\begin{align}
M\int^T_0\int_{\Sigma}\phi^i_t(t,Ct,z)dz+M\int_{\Sigma}\phi^i(0,0,z)dz+M^2\sum^d_{k=1}\frac{\cosh(\frac{1}{2})}{\sinh(\frac{1}{2})}\int^T_0\int_{\Sigma}\partial_k\phi^i(t,Ct,z)dz=0.
\end{align}
From (\ref{5.5}) and (\ref{5.11}), we deduce that
\begin{multline}
\int^T_0\int_{\mathbb{T}^d}(u\cdot\phi_t+\nabla u:\nabla \phi_t) dxdt+\int_{\mathbb{T}^d}(u_0\cdot\phi(0,x)+\nabla u_0:\nabla \phi(0,x))dx\\+\int^T_0\int_{\mathbb{T}^d}T^a:\nabla \phi dxdt
+\sum^d_{i,j,k=1}\int^T_0\int_{\mathbb{T}^d}u^j\partial_ku^i\partial_{jk}\phi^i dxdt=0,
\end{multline}
which implies that (\ref{5.1}) is the weak solution of the system (\ref{1}).

\begin{rema}
In the whole space, the Euler-Poincar\'e equations possess a special class of peakon solutions as follows:
\begin{align}
u^i=Me^{-|a\cdot x-Ct|}=Me^{-|\sum^d_{j=1}a^jx^j-Ct|}, i=1,2,...d.
\end{align}
One can verify that the above solutions are weak solutions of the Euler-Poincar\'e equations in the sense of distributions by the above similar argument.
\end{rema}

 {\bf Acknowledgements}. This work was
partially supported by NNSFC (No.11671407 and No.11701586),  FDCT (No. 098/2013/A3), Guangdong Special Support Program (No. 8-2015), and the key project of NSF of  Guangdong province (No. 2016A030311004).
%The authors thank the referees for their valuable comments and suggestions.

\bibliographystyle{abbrv} %plain ,%alpha, %abbrv
\bibliography{myref}
\end{document}